%% file: Guy-Roudolf-MPRF-V2.tex
\documentclass[a4paper,11pt,twoside]{paper}
\usepackage[latin1]{inputenc}
\usepackage{amsmath,amssymb,amsthm,mathrsfs,,mathtools,graphicx}
\usepackage{hyperref}
\newtheoremstyle{sans}{\parskip}{\parskip}{\itshape}
                       {0pt}{\bfseries\sffamily}{.}{ }{}
\newtheoremstyle{sansplain}{\parskip}{\parskip}{}
                       {0pt}{\bfseries\sffamily}{.}{ }{}

\theoremstyle{sans}
\newtheorem{prop}{Proposition}[section]
\newtheorem{coro}[prop]{Corollary}
\newtheorem{thm}[prop]{Theorem}
\newtheorem{lem}[prop]{Lemma}

\theoremstyle{sansplain}
\newtheorem{rem}[prop]{Remark}
\newtheorem{defin}[prop]{Definition}

\makeatletter
\@addtoreset{equation}{section}
\makeatother

\newcommand\C{\mathbb{C}}
\newcommand\Pb{\mathbb{P}}
\newcommand\Sc{\mathcal{S}}
\newcommand\Zb{\mathbb{Z}}

\renewcommand{\geq}{\geqslant}
\renewcommand{\leq}{\leqslant}

\def\H{\mathcal{H}}

\def\DD{\displaystyle}
\def\egaldef{\stackrel{\mbox{\tiny def}}{=}}
\renewcommand\>[1]{\vec{#1}}
\newcommand\wt[1]{\widetilde{#1}}
\newcommand{\dproof}{\noindent {\it Proof.} \quad}
\newcommand{\fproof}{\hfill $\blacksquare$ }
\begin{document}
\title{Random Walks in the Quarter-Plane: Advances in  Explicit Criterions for the Finiteness
 of the Associated Group in the Genus 1 Case}
\author{Guy Fayolle\thanks{INRIA Paris-Rocquencourt, Domaine de Voluceau, BP 105, 78153 Le Chesnay Cedex, France. Email: {\tt Guy.Fayolle@inria.fr}}    \and
        Roudolf Iasnogorodski\thanks{Saint-Petersbourg, Russia.
       Email: \texttt{iasnogorodski@mail.ru}}}
\date{\today}
\maketitle
\begin{abstract}
In  the book \cite{FIM}, original methods were proposed to determine the invariant measure of random walks in the quarter plane with small jumps, the general solution being obtained via reduction to boundary value problems. Among other things, an important quantity, the so-called \emph{group of the walk}, allows to deduce theoretical features about the nature of the solutions. In particular, when the \emph{order} of the group is finite, necessary and sufficient conditions have been given in \cite{FIM} for the solution to be rational or algebraic. In this paper,  when the underlying algebraic curve is of  genus $1$, we propose a concrete criterion ensuring the finiteness of the group. It turns out that this criterion can be expressed as the cancellation of a determinant of a matrix of order $3$ or $4$, which depends in a polynomial way on the coefficients of the walk. \end{abstract}


AMS $2000$ Subject Classification: Primary 60G50; secondary 30F10, 30D05

\section{Introduction}
In a probabilistic framework, we consider a piecewise homogeneous random walk with sample paths in $\Zb_+^2$, the lattice in the positive quarter plane. In the strict interior of $\Zb_+^2$,  the size of the jumps is $1$,  and  $\{p_{ij}, |i|,|j| \leq 1\}$ will denote the generator of the process for this region. Thus a transition $(m,n)\to(m+i,n+j), mn>0,$ can take place with probability $p_{ij}$, and
\[
\sum_{ |i|,|j| \leq 1} p_{ij} =1.
\]
On the other hand, no strong assumption is made about the boundedness of the upward jumps on the axes, neither at $(0,0)$. In addition, the downward jumps on the $x$ [resp.\ $y$] axis are bounded by $L$ [resp.\ $M$], where $L$ and $M$ are arbitrary finite integers. 

Needless to recall that a huge amount of work has been devoted to the analysis of this process since the early 1970s, both from analytic and probabilistic points of view. Today, it can be asserted that the main issues at stake (e.g. ergodicity conditions, computation of the invariant measure, etc) have been settled.  

In the book \cite{FIM}, an important quantity was studied, the so-called \emph{group of the walk}, originally introduced in \cite{MALY70,MALY71c}. When the \emph{order} of this group is finite, the nature of the solutions can be fully characterized. In particular, necessary and sufficient conditions have been obtained in \cite{FIM} for these solutions to be rational or algebraic. In brief, as explained hereafter in Sections (\ref{EQUA}) and (\ref{GENUS}), this group exchanges the two roots of the algebraic curve defined by $Q(x,y)=0$ (see equation \eqref{eq:eqfonc}), the  \emph{genus} of which is either $0$ or $1$. 

 In a combinatorial context (enumeration of lattice walks), it was possible, by following the approach of \cite{FIM}, to analyze the nature of the bivariate counting generating functions and to discriminate between algebraicity and holonomy.  

On the other hand, when the genus of the random walk  is equal to $0$, an effective criterion giving the order of the group has been provided in \cite{FR2}: then the group is infinite, except  precisely when the drift vector of the walk in the interior of $\Zb_+^2$ is equal to $0$, where finiteness is quite possible.

The genus $1$ case is more difficult and  solved in this paper, in the sense that we propose explicit necessary and sufficient criteria for the finiteness of the group. 

\subsection{The basic functional equations (see \cite{FIM}, chapters 2 and 5)}\label{EQUA}
 The invariant measure $\{\pi_{i,j}, i,j \ge 0\}$  does satisfy the fundamental functional equation 
\begin{equation} \label{eq:eqfonc}
 Q(x,y) \pi(x,y) = q(x,y) \pi(x) + \wt{q}(x,y)\wt{\pi}(y) + \pi_0(x,y), 
\end{equation}
with
\begin{equation*}  
\label{A2-FE}
\begin{cases}
\pi(x,y) = \DD \sum_{i,j \geq 1} \pi_{ij} x^{i-1} y^{j-1}, \\[0.5cm]  
\pi(x) = \DD \sum_{i\geq L} \pi_{i0} x^{i-L}, \quad \wt{\pi}(y) = \sum_{j\geq M} \pi_{0j} y^{j-M}, \\[0.5cm] 
Q(x,y) = \DD xy \Bigg[ 1 - \sum_{i,j\in\Sc} p_{ij} x^i y^j \Bigg], \quad  \sum_{i,j\in\Sc} p_{ij} =1, \\[0.5cm] 
q(x,y) = \DD  x^L \Bigg[\sum_{i \geq -L, j \geq 0}
p'_{ij} x^i y^j - 1 \Bigg] \equiv x^L(P_{L0}(x,y)-1), \\[0.5cm]
\wt{q}(x,y) = \DD y^M \Bigg[\sum_{i \geq 0, j \geq -M} p''_{ij} x^i y^j - 1 \Bigg] 
\equiv y^M(P_{0M}(x,y)-1), \\[0.6cm]
 \pi_0(x,y) = \DD  \sum_{i=1} ^{ L-1}\pi_{i0} x^i\big[P_{i0}(x,y)-1\big] 
 +  \sum_{j=1} ^{ M-1}\pi_{0j} y^j\big[P_{0j}(x,y)-1\big] + \pi_{00}(P_{00}(xy)-1).
  \end{cases}
\end{equation*}
In equation \eqref{eq:eqfonc}, $\mathcal{S}$ is the set of allowed jumps, the unknown functions 
$\pi(x,y), \pi(x), \wt{\pi}(y)$ are sought to be analytic in the region $\{(x,y)\in \mathbb{C}^{2} : |x|<1,|y|<1\}$, and continuous on their respective boundaries. In addition, $q, \wt{q}, q_0, P_{i0}, P_{0j}, $ are given probability generating functions supposed to have suitable analytic continuations (as a rule, they are polynomials when the jumps are bounded). 

The function  $Q(x,y)$, often referred to as the \emph{kernel} of (\ref{eq:eqfonc}), can be rewritten in the two following equivalent forms
     \begin{equation}
     \label{eq:kernel}
          Q(x,y) = a(x) y^{2}+ b(x) y + c(x) = \wt{a}(y) x^{2}+
          \wt{b}(y) x + \wt{c}(y),
     \end{equation}
where
\begin{alignat*}{2} \label{eq:abc}
a(x) &=p_{1,1}x^{2}+p_{0,1}x+p_{-1,1}  &  \qquad \wt{a}(y) & = p_{1,1}y^{2}+p_{1,0}y+p_{1,-1}, \\
b(x) &= p_{1,0}x^{2}+(p_{0,0}-1)x+p_{-1,0}  & \qquad \wt{b}(y) & = p_{0,1}y^{2}+(p_{0,0}-1)y +p_{0,-1}, \\
c(x) &=p_{1,-1}x^{2}+p_{0,-1}x+p_{-1,-1}  &  \qquad \wt{c}(y) & = p_{-1,1}y^{2}+p_{-1,0}y+p_{-1,-1}.
\end{alignat*}
We shall also need the discriminants
     \begin{equation}\label{eq:discri}
    D(x) \egaldef b^{2}(x)-4a(x) c(x),  \qquad  \wt{D}(y) \egaldef \wt{b}^{2}(y)-4 \wt{a}(y)\wt{c}(y).
     \end{equation}     
The polynomials $D$ and $\wt{D}$  are of degree $4$, respectively in $x$ and $y$, with dominant coefficients
     \begin{equation} \label{eq_C}
          d_4 = p_{1,0}^{2}-4p_{1,1}p_{1,-1},
          \qquad \wt{d_4}=p_{0,1}^{2}-4p_{1,1}p_{-1,1}.
     \end{equation}
More information is given in Appendix \ref{app:unif}
 \subsection{The group and the genus}\label{GENUS} 
 
Let $\C(x), \C(y)$ and $\C\,(x,y)$ denote the respective fields of rational functions of $x, y$ and $(x,y)$ over $\C$. Since in general $Q$ is assumed to be irreducible, the quotient field $\C\,(x,y)$ with respect to $Q$ is also a field and will be denoted by $\C_Q(x,y)$.
 \begin{defin}
The \emph{group of the random walk} is the Galois group  $\H=\langle \xi,\eta\rangle$ of automorphisms of   $\C_Q(x,y)$ generated by  $\xi$ and $\eta$ given by
     \begin{equation*}
          \xi(x,y)=\Bigg(x,\frac{c(x)}{y\,a(x)}\Bigg),\qquad 
          \eta(x,y)=\Bigg(\frac{\wt{c}(y)}{x\,\wt{a}(y)},y\Bigg).
     \end{equation*}
\end{defin}
Here, $\xi$ and $\eta$ are involutions, satisfying  $\xi^2=\eta^2=I$.

 Let
\begin{equation*}
\delta =\xi\eta,
\end{equation*}
the non commutative product. Then $\H$ has a normal cyclic subgroup $\H_0=\{\delta^i, i\in\Zb\,\}$, which is finite or infinite, and 
$\H / \H_0$ is a group of order 2. 

Let $I$ denote the identity operator. Then the group $\H$ is finite of order $2n$ if, and only if, 
\begin{equation}\label{eq:gr2n}
 \delta^n = I.
\end{equation}
We shall write $f_\alpha=\alpha(f)$ for any automorphism $\alpha\in\H$ and any function 
$f\in\C_Q(x,y)$.  As briefly depicted in Section \ref{app:unif} of the Appendix, the fundamental equation \eqref{eq:eqfonc}, together with $\xi,\eta,\delta$, can be \emph{lifted onto the universal covering} $\C$ (the finite complex plane).

In \cite{FIM}, the group $\H$ of the random walk was shown
to be of even order $2n, n= 2,\ldots,\infty$, and we were able to characterize completely the solutions of the basic functional equation. Moreover, when  $n$ is supposed to be finite and the functions
$q,\wt{q},q_0$ are polynomials, we also gave necessary and sufficient conditions for these solutions to be rational or algebraic.

\section{Topics on the conditions for $\H$ to be finite \label{sec:finite}}
 Finding an exact explicit form for  $n$ to be finite is a deep question, which has many connections with some classical problems in algebraic geometry. As quoted in equation \eqref{eq:gfinite} of Appendix\ref{app:unif}, a necessary and sufficient condition was provided in \cite{FIM}, which, albeit  theoretically nice (!) is not really easy to check by calculus.The main goal of the present paper is precisely to transform this condition into closed-form expressions, algebraically tractable. As we shall see, there are some structural differences due to the parity of the number  $n$ defined in Section \ref{GENUS}, remembering that the order of the group is $2n$. Therefore, as a kind of introductory purpose , we  will consider first the groups of order $4, 6, 8$.

 Recalling that $\H$ is generated by the elements $\xi$ and
$\eta$,  we can define  the homomorphism 
$$ h(R(x,y)) \egaldef R(h(x),h(y)), \quad \forall h\in\H, \ \forall
R\in\C\,_Q(x,y).$$ Clearly, $2$ elements $h_1,h_2$ of $\H$ are
identical if, and only if,
$$ h_1(x) = h_2(x), \quad h_1(y) = h_2(y) .$$ In addition, for any $ R
\in\C\,_Q(x,y)$, the following important equivalences hold:
\begin{equation} \label{eq4.1.0}\begin{cases}
\xi(R) = R \Longleftrightarrow \quad R\in\C\,(x) ,\\ \eta(R) = R
\Longleftrightarrow \quad R\in\C\,(y) ,
\end{cases} 
\end{equation}
so that $\C\,(x)$ (resp. $\C\,(y)$ are the elements of $\C\,_Q(x,y)$
invariant with respect to $\xi$ (resp. $\eta$). Indeed, let

\begin{equation} \label{eq4.1.00} 
\rho (x,y) = \frac{P_0(x,y)}{P_1(x,y)}, 
\end{equation} 
where $P_0(x,y),P_0(x,y)$ are polynomials in $x,y$. Since 
 $Q(x,y)$ is a polynomial of second degree in $y$, 
\eqref{eq4.1.00} yields at once
\[
\rho (x,y) = \dfrac{A_1(x)y + A_0(x)}{B_1(x)y + B_0(x)} \quad \pmod{Q(x,y)}.
\]
When $\rho$ is supposed to be invariant with respect to $\xi$,
it follows that
\begin{equation} \label{eq4.1.01} \dfrac{A_1(x)y + A_0(x)}{B_1(x)y +
B_0(x)} = \frac{A_1(x)y_{\xi} + A_0(x)}{B_1(x)y_{\xi} + B_0(x)} \quad  \pmod{Q(x,y)}. 
 \end{equation} 
 As $y\neq y_{\xi}$, we obtain from  (\ref{eq4.1.01})
$$A_1(x) B_0(x) = A_0(x) B_1(x), $$ 
so that
\begin{equation} \label{eq4.1.02}
   \rho (x,y)\equiv\rho (x) = \begin{cases} \dfrac{A_0(x)}{B_0(x)}, & \quad \text{if}
\quad B_0 \neq 0 , \\[0.3cm] \dfrac{A_1(x)}{B_1(x)} & \quad
\text{otherwise},
\end{cases}
 \end{equation} 
 which proves (\ref{eq4.1.0}).

\begin{defin} Introduce the matrix 
\begin{equation}\label{eq:mat} 
\Pb= \begin{pmatrix}
p_{11} & p_{10} & p_{1,-1} \\ 
p_{01} & p_{00}-1 & p_{0,-1} \\
p_{-1,1} & p_{-1,0} & p_{-1,-1} 
\end{pmatrix}, 
\end{equation}
and let $\>C_1, \>C_2, \>C_3$  (resp.  $\>D_1, \>D_2, \>D_3$) denote the column vectors of $P$ (resp. of $\Pb^T$, the tranpose matrix of $\Pb$).
\fproof
\end{defin} 

The simple following property will  be very useful.
\begin{prop} \label{order2s}
Assume there exists a positive integer $s$, such that
\begin{equation} \label{eq:gr2s}
\delta^s(x)=x.
\end{equation}
Then $\delta^s=I$ and the group is of order $2s$, where  $s$ stands for the smallest integer with property 
\eqref{eq:gr2s}.
\end{prop}
\dproof
Each of the three following  permutations
\[
x\Longleftrightarrow y, \quad \delta \Longleftrightarrow \delta^{-1}, \quad \Pb \Longleftrightarrow \Pb^T,
\]
 implies the two other ones.  Hence, the quantity $\rho(x,y,k)\egaldef\delta^k(x). \delta^{-k}(y)$, for any integer $k\ge1$, remains invariant by permuting $\Pb$ with $\Pb^T$. 

 Assume first $s=2m$. Then \eqref{eq:gr2s} becomes  $\delta^m(x)=\delta^{-m}(x)$, and
\[
 \rho(x,y,m) = \delta^{-m}(x). \delta^{-m}(y) =  \delta^m(x). \delta^m(y),
 \]
where the second equality is obtained  by replacing $\Pb$ by $\Pb^T$. Then, comparing with the definition of $\rho(x,y,m)$, we get $\delta^m(y)=\delta^{-m}(y)$, which yields in turn $\delta^s(y)=y$, whence $\delta^s=I$.

The argument works exactly the same way if $s$ is odd, say $s=2m+1$. Indeed, in this case we have 
\[
 \rho(x,y,m) = \delta^{-(m+1)}(x). \delta^{-m}(y) =  \delta^m(x). \delta^{m+1}(y),
 \]
(by exchanging again  $\Pb$ with $\Pb^T$), which implies $\delta^{m+1}(y)= \delta^{-m}(y)$, that is $\delta^s(y)= y$, concluding the proof of the proposition. \fproof
\medskip
\begin{coro}\label{orders} \mbox{ }
\begin{enumerate}
\item \label{pt1} If there exists $s$ such that $\delta^s(x)=r(x)$, where $r(x)$ represents a rational fraction of 
$x$, then $\delta^{2s}(x)=x$ and the group is of order $4s$.
\item \label{pt2} If there exists $s$ such that $\delta^s(x)=t(y)$, where $t(y)$ represents a rational fraction of 
$y$, then $\delta^{2s+1}(x)=x$ and the group is of order $4s+2$.
\end{enumerate}
\end{coro}
\dproof
Remark first the identities $\xi \delta^s\xi= \delta^{-s}$ and $\eta\delta^s\xi= \delta^{-s-1}$. 

So, we have the following chain of equalities. 
\[
\delta^s(x) = r(x) \Longrightarrow \xi \delta^s\xi(x)= \delta^s(x)\Longleftrightarrow  \delta^{-s}(x)=\delta^s(x) \Longleftrightarrow  \delta^{2s}(x) = x.
\]
Similarly
\[
\delta^s(x) = t(y) \Longrightarrow \eta\delta^s\xi(x)= \delta^s(x)\Longleftrightarrow  \delta^{-s-1}(x)=\delta^s(x) \Longleftrightarrow  \delta^{2s+1}(x) = x.
\]
In both cases, the conclusion follow from Proposition \ref{order2s}, in which $s$ is replaced respectively by $2s$ and $2s+1$.
\fproof
\medskip
\begin{lem} \label{lem:frac1}
On the algebraic curve $\{Q(x,y) = 0\}$, the following general relations hold.
\begin{equation} \label{eq4.1.6}
\left\{ \begin{array}{lll} 
\eta(x) &=& {\displaystyle  \frac{xv(y) - u(y)}{xw(y) - v(y)},} \\
\ \\ \xi(y) &=& {\displaystyle  \frac{y \wt{v}(x) - \wt{u}(x)}{y \wt{w}(x) - \wt{v}(x)} ,} 
\end{array} \right.
\end{equation} 
where $u, v, w, h$ (resp. $ \wt{v}, \wt{v}, \wt{w}, \wt{h}$) are polynomials of degree $\le2$. 

In particular, there exist affine solutions
\begin{equation}\label{eq:affine}
(u(y),v(y),w(y))^T = \>A\,y+ \>B, \  (\wt{u}(x),\wt{v}(x),\wt{w}(x))^T =  \>E\,x + \>F,
\end{equation}
with  column vectors
\begin{equation*}\label{eq:vec}
\>A = (u_0,v_0,w_0)^T, \ \>B = (u_1,v_1,w_1)^T, \ 
\>E = (\wt{u}_0,\wt{v}_0,\wt{w}_0)^T, \ 
\>F = (\wt{u}_1, \wt{v}_1,\wt{w}_1)^T,
\end{equation*}
given by
\begin{equation}\label{eq:sys}
\left\{ \begin{array}{lll}
\>A  & = & (\alpha\>C_2 + \beta \>C_1)\times \>C_3, \\[0.1cm]
\>B  & = & \>C_1\times(\alpha\>C_3 +\beta\>C_2), \\[0.1cm]
\>E  & = &  (\wt{\alpha}\>D_2 + \wt{\beta} \>D_1)\times \>D_3, \\[0.1cm]
\>F  & = &  \>D_1\times(\wt{\alpha}\>D_3 + \wt{\beta} \>D_2), 
\end{array}\right.
\end{equation}
where $\alpha, \wt{\alpha}, \beta, \wt{\beta}$ are arbitrary complex constants,
and the operator ``$\,\times$'' stands for the cross vector product. In addition, when $\Pb$ is of rank $3$, none of the vectors $\>A, \>B, \>E, \>F$ do not vanish.  Choosing in  \eqref{eq:sys}   \mbox{$\alpha=\wt{\alpha}=0, \beta=\wt{\beta}=1$}, gives

\begin{equation}\label{eq:special}
\left\{ \begin{array}{l}
u(y)  =  y \Delta_{13} -\Delta_{12},  \qquad \wt{u}(x)  =  x \Delta_{31}-\Delta_{21}, \\[0.1cm]
v(y)  =  y \Delta_{23} -\Delta_{22},  \qquad \wt{v}(x)  =  x \Delta_{32}-\Delta_{22}, \\[0.1cm] 
w(y) =  y \Delta_{33}-\Delta_{32},  \qquad \wt{w}(x)  =   x \Delta_{33}-\Delta_{23},
\end{array}\right.
\end{equation}
 where $\Delta_{ij}$ denotes the cofactor of the $(i,j)^{th}$ entry of the matrix 
 $\Pb$ given in \eqref{eq:mat}.
\end{lem} 

\dproof
We proceed by construction, assuming $\eta(x)$ is given by the following expression
\begin{equation}\label{eq:etainv}
\eta(x) =  \frac{xv(y) - u(y)}{xw(y) - h(y)},
\end{equation}
where $u, v, w, h$ are affine functions of $y$. Then, by using the basic relation
\[
\eta(x) = \frac{\wt{c}(y)}{x \wt{a}(y)}
\]
 in \eqref{eq:etainv}, we obtain 
\[
\wt{a}(y)v(y) x^2 - [\wt{a}(y)u(y) + \wt{c}(y)w(y) ]x +\wt{c}(y)h(y) =0 \ \pmod{Q(x,y)},
\]
which must  be proportional to \eqref{eq:kernel} written in the form
\[
\wt{a}(y) x^2 + \wt{b}(y)x +\wt{c}(y) =0.
\]
This implies the two identities,  $\forall y \in \C$,
\begin{equation}\label{eq:autom1}
\begin{cases}
\wt{a}(y)u(y) + \wt{b}(y)v(y) +\wt{c}(y)w(y) = 0, \\[0.1cm]
h(y)=v(y).
\end{cases}
\end{equation}

Hereafter, we assume the matrix $\Pb$ is of rank $3$. Indeed, it will be shown in the next section that   \emph{this is always the case, except when the group is of order $4$}]. 

So,  letting $\mathcal{V}$ denote the vector space of polynomials of degree $\le2$, the polynoms $\wt{a},\wt{b},\wt{c}$ form a base of $\mathcal{V}$, and we can write
\begin{equation*}
\begin{pmatrix}
u\\v\\w
\end{pmatrix} =
\mathcal{M}
\begin{pmatrix}
\wt{a}\\ \wt{b} \\ \wt{c})
\end{pmatrix}, 
\end{equation*}
where $\mathcal{M}$ is an unspecified constant  matrix. Then the first equation of \eqref{eq:autom1} says that one must look for elements $(u, v, w)\in\mathcal{V}^3$, such that the two vectors  
$(u, v, w)$ and $(\wt{a},\wt{b},\wt{c})$ are \emph{orthogonal}. The bilinear mapping $(u, v, w).(\wt{a},\wt{b},\wt{c})^T$ gives rise to an associated quadratic form (in the variables $\wt{a},\wt{b},\wt{c}$) 
\[
B(\wt{a},\wt{b},\wt{c})=(\wt{a},\wt{b},\wt{c})\, \mathcal{M}^*(\wt{a},\wt{b},\wt{c})^T,
\]
where $\mathcal{M}^*$ is a symmetric matrix. Consequently, since $\wt{a}(y),\wt{b}(y),\wt{c}(y)$ build a base of $\mathcal{V}^3$,  $B(.)$ will be identically zero, $\forall y\in\C$, if, and only if, 
$\mathcal{M}^*=0$, which implies that $\mathcal{M}$ is a skew-symmetric matrix
\[
\mathcal{M}=
 \begin{pmatrix}
0 & -\varepsilon_3 & \varepsilon_2 \\ 
\varepsilon_3 & 0 & -\varepsilon_1 \\
-\varepsilon_2 & \varepsilon_1 & 0
 \end{pmatrix}.
 \]
Next, we  can characterize the $3$  basic families $u(y),v(y),w(y)$, which among other things generate \eqref{eq:affine} by linear combination.
\begin{enumerate}
\item[1.] $u,v,w$ of degree $1$.  It suffices to choose 
\[
\varepsilon_3 = p_{-1,1}, \ \varepsilon_2 = p_{0,1}, \ \varepsilon_1 =  p_{1,1}.
\]
\item[2.] $u,v,w$ without constant terms, so that $u/y,v/y,w/y$ are admissible and still of degree~$1$. 
It suffices to choose 
\[
\varepsilon_3 = p_{-1,-1}, \ \varepsilon_2 = p_{0,-1}, \ \varepsilon_1 =  p_{1,-1}.
\]
\item[3.] $u,v,w$ of degree $2$. It suffices to choose 
\[
\varepsilon_3 = p_{-1,0}, \ \varepsilon_2 = p_{0,0}-1, \ \varepsilon_1 =  p_{1,0}.
\]
\end{enumerate}
The proof of the lemma is terminated.
\fproof
\medskip
\begin{lem} \label{lem:frac2}
Let $\gamma$ an endomorphism defined on the algebraic surface $\{Q(x,y)=0\}$, which is supposed to be invariant on the field $\C(x)$ of rational functions of $x$, and such that 
\begin{equation}\label{eq:xi1}
\gamma(y) = \frac{yf(x)-e(x)}{yg(x)+h(x)},
\end{equation}
where $e,f,g,h$ are polynomials of degree $1$ in $x$. (Note that this is always possible, as shown in 
Lemma~\ref{lem:frac1}). Then, for $\gamma$ to be an involution, the condition \mbox{$f(x)+h(x)\equiv0$} is necessary and sufficient.
\end{lem}
\dproof Applying $\gamma$ to both terms of \eqref{eq:xi1} yields the equality
\[
\gamma^2(y) = \frac{y[f^2(x)-e(x)g(x)] - e(x)[f(x)+h(x)]}{yg(x)[f(x)+h(x)] - e(x)g(x)+h^2(x)}.
\]
If $f(x)+h(x)=0$, then we have immediately $\gamma^2(y)=y$, and consequently $\gamma^2=I$, showing at once that $f(x)+h(x)=0$ is a sufficient condition for $\gamma$ to be an involution.

On the other hand, we have
\[
g(x)(f(x)+h(x))y^2 + (h^2(x)-f^2(x))y + e(x)(f(x)+h(x))=0.
\]
Comparing the last equation with \eqref{eq:kernel}, we obtain  (omitting the variable $x$)
\begin{equation}\label{eq:ratio}
\frac{g(f+h)}{a} = \frac{h^2-f^2}{b}= \frac{e(f+h)}{c}, \  \pmod{Q(x,y)}.
\end{equation}
 Assume for a while $f+h\not\equiv0$. Then, by \eqref{eq:ratio},
\begin{equation}\label{eq:ratio1}
\begin{cases}
bg=a(h-f),\\
cg=ae.
\end{cases}
\end{equation}
The second degree polynomials $a(x)$ and $b(x)$ are relatively prime. Indeed, the roots of $a(x)$ are either both negative, or complex conjugate. But $b(x)$ does not admit roots with a negative real part, since $p_{00}-1 <0$. So, $a(x)$ should divide $g(x)$, which is impossible. 

So, in equation \eqref{eq:xi1}, we must have $f(x)+h(x)=0, \forall x\in\C$, except when $a,b,c$ are all of degree $1$, which corresponds to the \emph{singular random walk} $p_{11}=p_{10}=p_{1,-1}=0$ introduced in \cite{FIM}. But we note in this latter case, that $Q(x,y)$ is of degree $1$ in $x$, and the genus of the algebraic curve is zero.   The lemma is proved. \fproof
\begin{rem}
The result of Lemma \ref{lem:frac2} does not hold if polynomials $e, f, h$ are not of degree $1$. For instance, one can check directly from \eqref{eq:ratio1} that, if $g$ or $e$ are taken to be of degree 2, then any involution $\gamma$ has necessarily  the form 
\[
\gamma(y)= \frac{(h-b) y - c}{a y + h},
\]
where $h$ is  an arbitrary polynomial, so that the condition $2h=b$ in general does not hold. 
\end{rem}

\subsection{Criterion for groups  of Order $4$ } \label{sec:gr4}

\begin{prop} \label{order4}$\quad$ The group $\H$ is of order 4 if, and only if,
\begin{equation}\label{eq:group4} 
 \begin{vmatrix}
p_{11} & p_{10} & p_{1,-1} \\ 
p_{01} & p_{00}-1 & p_{0,-1} \\
 p_{-1,1} & p_{-1,0} & p_{-1,-1} 
 \end{vmatrix} = 0, 
\end{equation}
and this is the only case where the matrix $\Pb$ has rank $2$.
\end{prop} 
\dproof Recalling that $\delta \egaldef \xi\eta$, the equality $\delta^2 = I$ can be rewritten as
$$\xi\eta = \eta\xi, $$ which by Proposition \ref{order2s} is for instance equivalent to 
\[
\xi\eta (x) = \eta(x),
\]
where we have used $\xi(x)=x$. So, $\eta(x)$ is left invariant by $\xi$, which implies
\[
\eta(x)\in \C\,(x).
\] 
Finally, $\eta$ is an involution ($\eta^2 = I $) and a conformal automorphism on $\C(x)$. Consequently,  $\eta$ is indeed a fractional linear transform of the type
\[
 \eta(x) = \dfrac{rx + s}{tx - r},
\]
 where all
coefficients belong to $\C\,$.  The following chain of equivalences hold.
\begin{eqnarray*}
\eta(x) = \dfrac{rx + s}{tx - r} &\Leftrightarrow & tx \eta(x) = r(x +
\eta(x)) + s \\ & \Leftrightarrow & 1, \; x + \eta(x), \; x \eta(x) \;
\mbox{are linearly dependent on $\C$} \\ & \Leftrightarrow & 1, \; -
\dfrac{\wt{b}(y)}{\wt{a}(y)}, \;
\dfrac{\wt{c}(y)}{\wt{a}(y)} \; \mbox{are linearly
dependent on $\C$} \\ & \Leftrightarrow & \wt{a}(y),
\wt{b}(y), \wt{c}(y) \; \mbox{are also linearly
dependent on $\C\,$,} \end{eqnarray*} where equation
\eqref{eq:kernel} has been used in the form
\[Q(x,y) = \wt{a}(y)x^2 + \wt{b}(y)x +
\wt{c}(y).
\] 
It is worth remarking that, starting from $\xi(y)$, the same argument would involve the transpose matrix $\Pb^T$, leading thus (as expected\,!) to the same criterion \eqref{eq:group4}.
 The proof of the lemma  is concluded. 
 \fproof

\subsection{Criterion  for groups  of Order $6$ } \label{sec:gr6}
\begin{prop} \label{order6}$\quad$ $\H$ is of order $6$ if, and only if,
\begin{equation} \label{eq:group6} 
\begin{vmatrix} 
 \Delta_{11} & \Delta_{21} & \Delta_{12} & \Delta_{22} \\
 \Delta_{12} & \Delta_{22} & \Delta_{13} & \Delta_{23} \\ 
 \Delta_{21} & \Delta_{31} & \Delta_{22} & \Delta_{32} \\
 \Delta_{22} & \Delta_{32} & \Delta_{23} & \Delta_{33}
\end{vmatrix} =0.
\end{equation}
 where the $\Delta_{ij}$'s have be given in Lemma \ref{lem:frac1}. 
 \end{prop}

 \dproof In this case $(\xi
\eta)^3 = I_d$, which is equivalent to
\begin{equation} \label{eq4.1.3} \eta \xi \eta = \xi \eta \xi. 
\end{equation}
Applying (\ref{eq4.1.3}) for instance to $x$, we get
$$\xi \eta(x) = \eta \xi \eta(x),$$ which shows that $\xi \eta(x)$ is
invariant with respect to $\eta$ and consequently is a rational
function of $y$, remembering one is dealing with the field of rational
functions. Similarly, $\eta \xi(y)$, invariant with respect to $\xi$,
is a rational function of $x$. Hence (\ref{eq4.1.3}) is plainly
equivalent to
\begin{align*} 
\begin{cases}
\xi \eta (x) = P(y), \\[0.1cm] \eta \xi(y) = R(x), 
\end{cases} 
\end{align*}
where $P$ and $R$ are rational. Then
 $$y = R(\xi \eta (x)) = R\circ P(y), $$ or, equivalently,
\begin{equation} \label{eq4.1.4}
R\circ P = I,
\end{equation}
so that $P$ and $R$ are  fractional linear transforms. Thus
(\ref{eq4.1.4}) yields the relation
\begin{equation}\label{eq4.1.5}
\xi(y) = \frac{p \eta(x) + q}{r \eta(x) + s} , \end{equation}
which imposes a linear dependence on $\C\,$ between the four elements
 1, $\xi(y)$, $\eta(x)$, $\xi(y)\eta(x)$, with $4$ unknown constants (in fact $3$ by homogeneity). Our goal is to avoid the pitfall of entering tedious (and harmful\,!) computations. 

Starting from equation \eqref{eq4.1.6}, we choose $\eta(x)$ by means of \eqref{eq:special}, that is 
\begin{equation}\label{eq:eta1}
\eta(x) =\frac{y(x\Delta_{23} - \Delta_{13}) -  x\Delta_{22} + \Delta_{12}}
{y(x\Delta_{33} - \Delta_{23}) -  x\Delta_{32} + \Delta_{22}}.
\end{equation}
Instantiating now \eqref{eq:eta1} in \eqref{eq4.1.5}, we obtain
\begin{equation}\label{eq:xi3}
\xi(y)= \frac{y[p(x\Delta_{23} - \Delta_{13})+q(x\Delta_{33} - \Delta_{23})] 
+ p(\Delta_{12}-x\Delta_{22})+q(\Delta_{22} - x\Delta_{32})}
{y[r(x\Delta_{23} - \Delta_{13}) +s(x\Delta_{33} - \Delta_{23})] 
+ r(\Delta_{12}-x\Delta_{22}) +s(\Delta_{22} - x\Delta_{32})}
\end{equation}
Then, rewriting \eqref{eq:xi3} in the  form proposed in \eqref{eq:xi1}, namely
\[
\xi(y) = \frac{yf(x)-e(x)}{yg(x)+h(x)},
\]
where
\begin{equation}\label{eq:xi4}
\left\{ \begin{array}{lll}
e(x)  & = & p(x\Delta_{22}-\Delta_{12}) + q(x\Delta_{32} -\Delta_{22}), \\
f(x)  & = &  p(x\Delta_{23} - \Delta_{13}) + q(x\Delta_{33} - \Delta_{23}), \\
g(x) & =  & r(x\Delta_{23} - \Delta_{13}) +s(x\Delta_{33} - \Delta_{23}), \\
h(x) & = &  r(\Delta_{12}-x\Delta_{22}) +s(\Delta_{22} - x\Delta_{32}),
\end{array}\right.
\end{equation}
we are in a position to compare system \eqref{eq:xi4} with the solution presented in  
equation~\eqref{eq:sys} of  Lemma~\ref{lem:frac1}. 

Indeed, letting $\>\varphi(x)= (e(x),f(x),g(x))^T\egaldef \>\varphi_0 + x \>\varphi_1$, we get from \eqref{eq:sys} the vector equation 
\[
\>\varphi_0+x\>\varphi_1 = \>E+x\>F,
\]
which, after renaming the constants $\wt\alpha$ and $\wt\beta$ as $c$ and $-d$ respectively, yields in turn $6$ linear equations, namely 
\begin{equation}\label{eq:xi5}
\begin{cases}
\ c\Delta_{11} + d\Delta_{21} + p\Delta_{12} + q\Delta_{22} = 0, \\
\ c\Delta_{12} + d\Delta_{22} + p\Delta_{13} + q\Delta_{23} = 0, \\
\ c\Delta_{13} + d\Delta_{23} + r\Delta_{13} + s\Delta_{23} = 0, \\
\ c\Delta_{21} + d\Delta_{31} + p\Delta_{22} + q\Delta_{32} = 0, \\
\ c\Delta_{22} + d\Delta_{32} + p\Delta_{23} + q\Delta_{33} = 0, \\
\ c\Delta_{23} + d\Delta_{33} + r\Delta_{23} + s\Delta_{33} = 0.
\end{cases}
\end{equation}

In addition, we have to take into account the constraint $f(x)+h(x)=0$, which is tantamount to 
\begin{equation}\label{eq:xi6}
 \begin{cases}
 p\Delta_{23} + q\Delta_{33} = r\Delta_{22} + s\Delta_{32}, \\[0.1cm]
 p\Delta_{13} + q\Delta_{23} = r\Delta_{12} + s\Delta_{22}.
\end{cases}
\end{equation}
So, the final step is to analyze the feasability of the global linear system formed by the intersection of \eqref{eq:xi5} and \eqref{eq:xi6}. Altogether, we are left with $8$ equations with respect to $6$ unknown variables $c,d,p,q,r,s$.

Equations of system \eqref{eq:xi5}, referred to as $1, 2\ldots 6$, can be split into two sets: 
\begin{itemize}
\item[(a)]  The set $(1,2,4,5)$, forming an homogeneous linear system of $4$ equations with 4 unknowns
\begin{equation}\label{eq:set4}
\begin{pmatrix}
\Delta_{11}&\Delta_{21}&\Delta_{12}&\Delta_{22}\\
\Delta_{12}&\Delta_{22}&\Delta_{13}&\Delta_{23}\\
\Delta_{21}&\Delta_{31}&\Delta_{22}&\Delta_{32}\\
\Delta_{22}&\Delta_{32}&\Delta_{23}&\Delta_{33}.
\end{pmatrix}
\begin{pmatrix}
c\\d\\p\\q
\end{pmatrix} =0
\end{equation}
\item[(b)] The set $(3,6)$, which can easily be rewritten as 
\begin{equation}\label{eq:set2}
 \begin{cases}
(c+r)\Delta_{13} + (d+s)\Delta_{23} = 0, \\[0.1cm]
(c+r)\Delta_{23} + (d+s)\Delta_{33} = 0.
\end{cases}
\end{equation}
\end{itemize}
Clearly, system \eqref{eq:set4} has a non trivial solution, if and only if condition \eqref{eq:group6} holds, which hence is necessary for the group to be of order $4$. To prove its sufficientcy, we have  to considers also systems \eqref{eq:set2}, \eqref{eq:xi6}.

The determinant of system \eqref{eq:set2} is equal to $\Delta_{13}\Delta_{33} -\Delta^2_{23}$. 
But the matrix $\Pb$, introduced in \eqref{eq:mat}, has all its entries $(i,j)$ positive, except $(2,2)=p_{00}-1<0$, so that
$$\Delta_{13}\Delta_{33} -\Delta^2_{23}\leq0,$$
Moreover, the equality $\Delta_{13}\Delta_{33} -\Delta^2_{23}=0$ takes place only for special values of the jump probabilities $p_{ij}$'s, corresponding to simple singular random walks (see \cite{FIM}), which we do not consider here. Consequently \eqref{eq:set2} has only the trivial solutions, so that
\begin{equation}\label{eq:set2bis}
c+r=d+s=0.
\end{equation}
Now, by \eqref{eq:set2bis}, one replaces $r$ and $s$ respectively by $-c$ and $-d$ in \eqref{eq:xi6}, and we get $2$ equations coïnciding in fact  with the equations $(2,5)$ of system \eqref{eq:xi5}. 
It is worth remarking that we found the two equations of \eqref{eq:xi6} are implicitly satisfied, but it was useful to check this fact, just for the sake of completeness\,! The proof of the lemma is concluded. 
\fproof

Now, we shall attack the general situation by slitting into the two possible situations $n=2m$ and $n=2m+1$. 
\subsection{Criterion for groups  of order $4m$ } \label{sec:gr4m}
In this rest of this section,  we refer without further comment to the notation and  formulae of Sections \ref{app:unif} and \ref{app:calculs} of  the Appendix. 
 
\begin{prop} \label{order4m} The group $\H$ is of order $4m$ if, and only if, the Weierstrass function 
$\wp$ with periods $(\omega_1,\omega_2)$ satisfies the equation 
\begin{equation}\label{eq:det-gr4n}
\wp(m\omega_3) = \wp(\omega_2/2).
\end{equation}
\end{prop}
\dproof 
Recalling that  $\delta=\xi\eta$, with $\xi(x)=x,\,\eta(y)=y,\,\xi^2=\eta^2=I$, 
we have here $\delta^{2m}=I$, 
that is 
\begin{equation}\label{eq:g2m-1}
(\xi\eta)^m = (\eta\xi)^m.
\end{equation} 
By applying equation \eqref{eq:g2m-1} at $x$ (or even at an arbitrary element of $\C(x)$), and replacing 
$x$ by $\xi(x)$ in the left-hand side member, we obtain
\[
\xi\delta^m(x) = \delta^m(x),
\]
which shows that the involution $\delta^m(x)$ is invariant with respect to $\xi$, and hence is an element of $\C(x)$. This can be summarized by
\begin{equation}\label{eq:g2m-2}
\delta^m(x) = F(x) =  \frac{xf-e}{xg-f},
\end{equation}
where $F(x)$ is a simple \emph{fractional linear transform}, with constants $e, f, g,$ to be determined.

From Proposition \ref{order2s} and Corollary \ref{orders}, it is worth recalling that  \eqref{eq:g2m-2} contains exactly the same information as condition \eqref{eq:g2m-1}. 

Thus, equation \eqref{eq:g2m-2} implies the existence of a linear dependence between the functions 
\begin{equation}\label{eq:indep}
x.\delta^m(x), \ x+\delta^m(x), \  \mathbf{1},
\end{equation}
where $\mathbf{1}$ denotes an arbitrary constant function and the symbol ``.''  dwells on the ordinary scalar product. Then, after a lifting onto the universal covering and a translation of $-m\omega_3/2$, condition \eqref{eq:indep} gives rise to the following lemma.

\begin{lem}\label{lem:indep}
For the group to be of order $4m$, a necessary and sufficient condition is that the  three functions 
 \begin{equation}\label{eq:indep1}
 x(\omega-m\omega_3/2).x(\omega+m\omega_3/2), \  
x(\omega-m\omega_3/2)+x(\omega+m\omega_3/2), \  \mathbf{1},
\end{equation}
be linearly dependent, $\forall\omega\in\C$. \fproof
\end{lem}
\dproof 
In agreement with \eqref{eq:SP}, let 
\begin{alignat*}{2}
S&=x(\omega+m\omega_3/2)+x(\omega-m\omega_3/2), 
\quad P=x(\omega+m\omega_3/2).x(\omega-m\omega_3/2), \\
X&=\wp(\omega), \quad Y=\wp(m\omega_3/2).
\end{alignat*}

By \eqref{eq:xunif},  $x(\omega)$ is homographic in $X$, and hence  Lemma \ref{lem:indep} amounts to saying that $S,P, \mathbf{1}$,  considered as functions of $X$, are linearly dependent. 

By system \eqref{eq:SP1} together with with definition \eqref{eq:SP1}, 
one verifies immediately that the existence of constants $e, f, g$ satisfying
\begin{equation}\label{eq:lin1}
eS +fP+g=0, \quad \forall X\in \C,
\end{equation}
 is synonymous with the  linear dependence
\begin{equation}\label{eq:lin2}
uA_1 + vB_1 + w(X-Y)^2 = 0, \quad \forall X\in \C,
\end{equation}
where the constants $u, v, w$ fulfil the linear sytem
\begin{equation}\label{eq:detuvw}
\begin{pmatrix}
u\\v\\w
\end{pmatrix} = 
\begin{pmatrix}
2p & p^2 & 1 \\
q-2rp & p(q-rp) & -r \\
 4r(pr-q) & 2(q-rp)^2 & 2r^2
\end{pmatrix}
\begin{pmatrix}
e \\ f \\ g
\end{pmatrix}.
\end{equation}
But the determinant of system \eqref{eq:detuvw} is exactly equal to $-2q^3$ 
(independent of $p$ and $r$, as expected) and does never vanish (see~\eqref{eq:xunif1}. Consequently, equations  
\eqref{eq:lin1} and \eqref{eq:lin2} are in fact equivalent. 

Ultimately, we have to extract the vector coefficients of $A_1,B_1,D$, which from Lemma \ref{lem:WSP} are \emph{polynomials of second degree} in $X$, so that condition \eqref{eq:indep1} amounts to 
\begin{equation}\label{eq:det-gr4n1}
\Delta(Y)\egaldef \ 
\begin{vmatrix}
4Y & 4Y^2-g_2 & -(g_2Y+2g_3) \\[0.1cm]
2Y^2 & g_2Y+2g_3 & 2g_3Y+g_2^2/8 \\[0.1cm]
1 & -2Y & Y^2
\end{vmatrix} = 0.
\end{equation}
The determinant \eqref{eq:det-gr4n1} is equal to
\[
8\,Y^6 -10\,g_2Y^4 - 40\,g_3Y^3 - 5/2\,g_2^2\,Y^2 - 2\,g_2\,g_3Y - 4g_3^2 + 1/8\,g_2^3,
\]
but it has the pleasant property of being the product of three explicit second degre polynomials, namely  
\begin{equation}\label{eq:poly}
\Delta(Y) = 8[Y^2 - 2e_1Y - (e_1^2 + e_2e_3)][Y^2 - 2e_2Y - (e_2^2 + e_3e_1)][Y^2 - 2e_3Y - (e_3^2 + e_1e_2)],
\end{equation}
where $e_1,e_2,e_3$ are defined in \eqref{eq:wp1} and \eqref{eq:wpe}.
Let $P_1,P_2,P_3$ denote the $3$ polynomials coming in \eqref{eq:poly}, respectively from left to right, with their corresponding reduced discriminants
\begin{alignat*}{2}
\delta_1 & = (e_1-e_2)(e_1-e_3)>0, \\
\delta_2 & = (e_2-e_3)(e_2-e_1)<0, \\ 
\delta_3 & = (e_3-e_1)(e_3-e_2)>0.
\end{alignat*}
 On the real interval $[0,\omega_2]$, $\wp(\omega)$ is positive and reaches its minimum at
  $\omega=\omega_2/2$, with  $\wp(\omega_2)=e_1$. So, we can immediately eliminate $P_2$, which has two complex roots. As for $P_3$, which has two real roots, one of them (at least) being negative, we have
  \[
  P_3(e_1) =  e_1^2 - 2e_1e_3 - e_3^2 -  e_1 e_2 = (e_1 -  e_2)( e_1 -  e_3) > 0.
\]
But $e_1$ is larger than the two roots of $P_3$, which therefore never cancel $\Delta(Y)$.
Thus we are left with the real roots in $Y$ of $P_1(Y)$, of which the sole one bigger than $e_1$ is admissible,  namely
\[
e_1+ \sqrt{(e_1 -  e_2)( e_1 -  e_3)} \equiv \wp\left(\frac{\omega_2}{4}\right).
\]
Finally, we have shown that \eqref{eq:det-gr4n1} holds if and only if 
\[
\wp(m\omega_3/2) = \wp(\omega_2/4).
\]
To conclude the proof of Proposition \ref{order4m}, it suffices to remark that, as we work inside the fundamental parallelogram $(\omega_1,\omega_2)$,  the  relations  
 $\frac{m\omega_3}{2} =  \pm\frac{\omega_2}{4}$ and $m\omega_3 =  \frac{\omega_2}{2}$ are equivalent 
 $\bmod\,(\omega_2)$.
 \fproof
\medskip

The calculation of $\wp(m\omega_3/2)$ could be carried out from \eqref{eq:swp}, via the recursive relationship
\[
 \wp((l+1)\omega_3/2) +  \wp((l-1)\omega_3/2) = 
\frac{(\wp(l\omega_3/2) + \wp(\omega_3/2))(4\wp(l\omega_3/2)\wp(\omega_3/2) -g_2)-2g_3}
{2(\wp(l\omega_3/2) - \wp(\omega_3/2))^2},
\]
siince the value of $\wp(\omega_3/2)$ is directly obtained from \eqref{eq:xunif} and \cite[Section3.3]{FIM}. Yet, one has to admit that the partial fraction giving $\wp(m\omega_3/2)$ in terms of   
$\wp(\omega_3/2)$ is hardly directly expoitable in Proposition  \ref{order4m}. 

Therefore, in the next sections, we shall pursue with the analysis of the condition \eqref{eq:indep}, having in mind the obtention of concrete formulas by means of  \emph{determinants} given in terms of the elements of the matrix $\Pb$ defined in \eqref{eq:mat}.
\subsubsection{The case $m=2k$}\label{sec:gr8k}
When $m=2k$, applying the operator $\delta^{-k}$ in \eqref{eq:indep} is equivalent to say, after a direct manipulation, that
 \begin{equation}\label{eq:g4m}
 \delta^k(x).\xi\delta^k(x), \ \delta^k(x) + \xi\delta^k(x), \  \mathbf{1},
 \end{equation}
are linearly dependent. But we know that $\delta^k(x).\xi\delta^k(x)$ and  $\delta^k(x)+\xi\delta^k(x)$ are elements of $\C_x$,  and by \eqref{eq:xunif1}, \eqref{eq:swp}, \eqref{eq:pwp}, they are in fact ratios of polynomials of degree 2 in~$x$.

In addition, letting $\zeta_j(x) \egaldef \delta^j(x)+\xi\delta^j(x)$, the following recursive scheme holds.
\begin{equation}\label{eq:rec-sx}
\begin{cases}
\zeta_0(x) & = 2x, \ \zeta_1(x) =  \eta(x) + \delta(x),\\[0.2cm]
\zeta_j(x) &= \zeta_{j-1}(\zeta_1(x))  - \zeta_{j-2}(x), \ \forall j\geq2.
\end{cases}
\end{equation}

\subsubsection{Explicit criterion for the groupe of order 8} \label{sec:gr8}
This correspondons to $m=2$  in the preceding section. In this case, explicit conditions can be carried out both manually and also with \emph{Maple~18}. Hereafter, we only list the main results.

In agreement with Section  \ref{sec:gr8k}, and using again $\xi^2=I$, the functions  
\[ 
\eta(x).\delta(x),  \  \eta(x) + \delta(x), \   \mathbf{1},
\]
are in $\C(x)$ and are sought to be linearly dependent on $\C$. Setting 
\[
 \eta(x) + \delta(x) = \frac{P_1(x)}{Q_1(x)}, \quad  \eta(x).\delta(x) = \frac{R_1(x)}{Q_1(x)},
 \]
the polynomials $P_1, Q_1, R_1$ satisfy the relations
\begin{equation} \label{eq:PQR8}
\begin{cases}
P_1  = -2\,S_2\,T_2 + S_3\,T_1 + S_1\,T_3, \\[0.1cm]
Q_1  = S_2^2 - S_3\, S_1, \\[0.1cm]
R_1  =  T_2^2 - T_3\,T_1,
\end{cases}
\end{equation}
with 
\begin{alignat*}{4}
S_1(x) & = p_{10}\,c(x) - p_{1,-1}\,b(x) & = & \ \  \Delta_{31}\,x - \Delta_{21} ,  \\
T_1(x) & = \frac{p_{-1,-1}\, b(x) - p_{-1,0}\,c(x)}{x} \ & = &  - \Delta_{21}\,x + \Delta_{11} ,\\[0.2cm]
S_2(x) & = p_{11} \,c(x)-p_{1,-1}\,a(x) & =  & - \Delta_{32}\,x + \Delta_{22} ,  \\
T_2(x) & = \frac{p_{11} \,a(x)-p_{-1,1}\,c(x)}{x} & = &\ \   \Delta_{22}\,x - \Delta_{12} ,\\[0.2cm]
S_3(x) & = p_{11}\,b(x) - p_{10}\,a(x) & = & \ \  \Delta_{33}\,x - \Delta_{23},   \\
T_3(x) & = \frac{p_{-1,0}\,a(x) - p_{-1,1}\,b(x)}{x} & = & - \Delta_{23}\,x + \Delta_{13} ,
\end{alignat*}
the $\Delta_{ij}$'s being the cofactors of the matrix  $\Pb$ introduced in Lemma \ref{lem:frac1}.
So, the condition ensuring the linear dependence between $P_1, Q_1, R_1$ given by \eqref{eq:PQR8}  leads to the following 
\smallskip

\begin{prop} \label{order8} $\quad$ The group $\H$ is of order $8$ if, and only if, the third order  determinant

\begin{equation} \label{eq:group8} 
 \begin{vmatrix} 
\begin{split} \MoveEqLeft[5] 
2\,\Delta_{22} \Delta_{32} \\ -(\Delta_{21} \Delta_{33} + \Delta_{31} \Delta_{23}) 
\end{split} 
  & \ \  \begin{split} 
 2\,(\Delta_{22}^2 - \Delta_{12} \Delta_{31} +\Delta_{21} \Delta_{23}) \\ 
 + \, \Delta_{11} \Delta_{33} + \Delta_{31} \Delta_{13} 
  \end{split} 
  & \begin{split} 
 2\,\Delta_{12} \Delta_{22}  \\ - (\Delta_{11} \Delta_{23} + \Delta_{21} \Delta_{13})  
 \end{split} \\[1cm]
  \Delta^2_{32} - \Delta_{31} \Delta_{33} 
  & -2\,\Delta_{32}\Delta_{22} +\Delta_{31} \Delta_{23} +\Delta_{21} \Delta_{33} 
  & \Delta^2_{22} - \Delta_{21} \Delta_{23} \\[0.1cm]
  \Delta^2_{22} - \Delta_{21} \Delta_{23} 
 & -2\,\Delta_{22}\Delta_{12} +\Delta_{11} \Delta_{23} +\Delta_{13} \Delta_{21}
 & \Delta^2_{12} - \Delta_{11} \Delta_{13} 
\end{vmatrix}
\end{equation}
is equal to zero. (Due to the size of the printing output, each  element of the first line of the matrix in \eqref{eq:group8} has been split vertically as the sum of two terms)  
 \fproof
 \end{prop} 

\medskip
 \begin{rem} It is interesting to note that the polynomials $S_i,T_i,i=1,2,3$, coïncide with  the ones appearing in \eqref{eq:special}.
\end{rem}
\bigskip
 \begin{rem} When $m=2$, i.e. a group of order $8$, the line of argument developed in Section \ref{sec:gr6} could be applied, but at the expense of  intricate computations. Indeed by Lemma \ref{lem:frac1},
\[
\eta(x) =  \frac{xv(y) - u(y)}{xw(y) - v(y)}, \qquad
 \xi\eta(x) = \frac{xv(\xi(y)) - u(\xi(y))}{xw(\xi(y)) - v(\xi(y))},
 \]
whence
\begin{alignat}{2}
\xi\eta(x)+\eta(x) & = \frac{(xv(\xi(y)) - u(\xi(y)))(xw(y) -v(y)) + (xv(y) -u(y))(xw(\xi(y)) - v(\xi(y)))}{(xw(\xi(y)) - v(\xi(y)))(xw(y)-v(y))}, \label{eq:s8}\\[0.2cm]
\xi\eta(x).\eta(x)  & =    \frac{(xv(\xi(y)) - u(\xi(y)))(xv(y) -u(y)}{(xw(\xi(y)) - v(\xi(y)))(xw(y)-v(y))}.  
\label{eq:p8}
  \end{alignat}
Setting for a while in  \eqref{eq:s8} and \eqref{eq:p8}
\[
 \xi\eta(x).\eta(x) = \frac{K_1}{L}, \quad  \xi\eta(x)+\eta(x) = \frac{K_2}{L},
\]
we must test a possible linear dependence between the functions
\[
K_1, \  K_2, \  L,
\]
which are all elements of $\C(x)$. Sketching the global calculus, we have for instance
\begin{alignat*}{2}
L & = (xw(y)-v(y))(xw(\xi(y)) - v(\xi(y)))\\
   & = w(y)w(\xi(y)) x^2 -x[v(y)w(\xi(y)) + w(y)v(\xi(y))] + v(y)v(\xi(y)).
   \end{alignat*}
A typical term coming in $L$ is  $ w(y)w(\xi(y))$, which from Lemma \ref{lem:frac1}  equals
 \begin{alignat*}{2}
 w(y)w(\xi(y)) & = w_0^2 +w_0w_1(y+\xi(y)) +w_1^2 y\xi(y)\\
 &= w_0^2 - w_0w_1 \frac{b(x)}{a(x)} + w_1^2 \, \frac{c(x)}{a(x)} . 
  \end{alignat*}
Analogously,   the other parts of the puzzle give 
\begin{alignat*}{2}
 a(x)L  =  & \ a(x)(xw_0-v_0)^2 - b(x)(xw_0-v_0)(xw_1-v_1) +c(x)(xw_1-v_1)^2, \\
 a(x)K_1  =  & \ a(x)(xv_0-u_0)^2 - b(x)(xv_0-u_0)(xv_1-u_1) +c(x)(xv_1-u_1)^2, \\
a(x)K_2 = & \  a(x)h_a(x) + b(x) h_b(x) + c(x)h_c(x),
\end{alignat*}
where
\begin{alignat*}{2}
h_a(x) = & \ 2 [v_0w_0x^2  - v_0^2x + v_0u_0], \\
h_b(x) = & \ -(v_1w_0-v_0w_1)x^2 +(2v_0v_1+u_0w_1+u_1w_0)x - (u_0v_1+u_1v_0),\\
h_c(x) = & \  2[v_1w_1x^2 -(v_1^2+u_1w_1)x + v_1u_1].
\end{alignat*}
At the very end,  this machinery would require to deal with rational fractions of \emph{4th~degree}, implying tedious computations, which are partly avoided by means of the procedure proposed above. \fproof
 \end{rem}
 
\subsection{Criterion for groups  of order $4m+2$} \label{sec:gr4m+2}
Here $m=2k+1$, and  condition \eqref{eq:gr2n} reads
\begin{equation}\label{eq:gr2k+1}
\delta^{2k+1}= I,
\end{equation}
which by Proposition \ref{order2s} and Corollary \ref{orders} is equivalent to
\[
\eta(\delta^k(x))  = \delta^k(x),
\]
or, by a simple algebra using  $\xi(x)=x$, to
\begin{equation}\label{eq:delta2k+1}
\delta^k(x)= G(y) \in\C(y).
\end{equation}
 Similarly, upon applying  \eqref{eq:gr2k+1} to $y$, we get
\[
 \delta^{-k}(y)=  \delta^{k+1}(y) =  \xi(\delta^{-k}(y),
\]
so that
\[
 \delta^{-k}(y) = F(x) \in\C(x).
\]
Then, applying  $\delta^{-k}$ to both members of  \eqref{eq:delta2k+1} yields
\[
x= \delta^{-k}(G(y)) = G( \delta^{-k}(y)) = G\circ F(x),
\]
which shows that $G\circ F=I$, and hence $G$ and $F$ are simple fractional linear transforms. 

Setting for instance
\[
G(y) = -\frac{py+q}{ry+s},
\]
 where $p,q,r,s$ are arbritrary complex constants, the problem is to achieve the linear relation
\begin{equation}\label{eq:gr4m+2}
r\,y\delta^k(x) + s\,\delta^k(x) + py +q=0  \  \pmod{Q(x,y)}.
\end{equation}

\section{Some general results}
We state below two concluding theorems allowing to conclude that, for the group to be finite,  there is a unique condition tantamout to the cancellation of a \emph{determinant}, the elements of which are intricate functions of the coefficients of the transition matrix $\Pb$, but nonetheless recursively computable. 

\begin{itemize}
\item The determinant is of order 3, for groups of order $4m, m\ge1$. 
\item The determinant is of order 4, for groups of order $4m+2,m\ge1$.
\end{itemize}

\subsection{A theorem about  $\delta^s$} \label{sec:general}
The following important fact holds. 

\begin{thm}\label{th:deltas}
For any integer $s\ge1$, we have
\begin{equation}\label{eq:deltas}
\delta^s(x)= \frac{y\,U_s(x) + V_s(x)}{W_s(x)} \  \pmod{Q(x,y)},
\end{equation}
where $U_s,V_s,W_s$ are second degree polynomials.
\end{thm} 
\dproof  For $s=1$, a direct (however slightly tedious) computation with the $(x,y)$ variables can be worked out. But, already for $s=2$, one seems to reach the limits of human computational abilities, which are definitely exceeded for $s\ge3$\,! A formal verification through the \emph{Maple~18} Computer Algebra System has been carried out for $s=2,3$, but seems hardly exploitable as soon as  $s>3$.  Hereafter, we propose a simple proof by using again the uniformisation \eqref{eq:xunif}.

Let
\[
\delta^s(x) - \xi\delta^s(x) \egaldef 2H, \ X= \wp(\omega), \ Y= \wp(s\omega_3).
\]
Then, with the notation of Appendix \ref{app:calculs},
\[
\delta^s(x)= \frac{S(\omega,s\omega_3)}{2} + H,
\]
where $S(\omega,s\omega_3)$ is given by \eqref{eq:SP1}.

Now, by \eqref{eq:swp}, \eqref{eq:pwp}, \eqref{eq:xunif1} and the addition formula for the $\wp$ function, we can write
\begin{alignat*}{2}
H& = \frac{q[\wp(s\omega_3- \omega) - \wp(s\omega_3 + \omega)]}
{2[\wp(s\omega_3 + \omega) -r][\wp(s\omega_3 -\omega)-r]} \\[0.2cm]
&= \frac{q\,\wp'(\omega)\wp'(s\omega_3)}
{2(X-Y)^2 [B(s\omega_3,\omega)-rA(s\omega_3,\omega) + r^2]} = 
\frac{q\,\wp'(\omega)\wp'(s\omega_3)}{D(X,Y)},
\end{alignat*}
 where $D(X,Y)$, introduced in \eqref{eq:SP2}, is a polynomial of second degree in $X$ and $Y$. On the other hand, by \eqref{eq:xunif} and \eqref{eq:xunif1}, we have 
\[
\wp'(\omega) = \frac{2q\,[2a(x)y+b(x)]}{(x-p)^2},
\]
and hence
\[
H= \frac{2q^2 \wp'(s\omega_3)[2a(x)y+b(x)]}{(x-p)^2D(X,Y)}.
\]
Since $x$ is a homographic function of $X$, it follows that the denominator of $H$ is a second degree polynomial in the variable $(x-p)$, where
\begin{alignat*}{2}
U_s(x) & = 4q^2\wp'(s\omega_3) a(x),  \\
V_s(x) & = 2q^2\wp'(s\omega_3) b(x) + (x-p)^2\left[2pB_1(X,Y)+ (q-2pr)A_1(X,Y) 
+ 4r(pr-q)(X-Y)^2\right],\\
W_s(x) & = (x-p)^2 D(X,Y).
\end{alignat*}
The proof of \eqref{eq:deltas} is concluded.
\fproof

\subsection{The form of the general criterion} 
\begin{thm}
For the group $\H$ to be finite, the  necessary and sufficient condition is 
$\det(\Omega)=0$, where $\Omega$ is a matrix of oder $3$ (resp.$4$) when the group is of order $4m$ (resp. 4m+2).
\end{thm}
\dproof We only present a sketch of the main line of argument.

\begin{itemize}
\item When the group $\H$  is of order $4m$, Theorem \ref{eq:deltas} and condition \eqref{eq:g4m} show that we are left with a system of $3$ homogeneous linear equations with $3$ unknowns.
\item  On the other hand, the situation for $\H$  to be of order $4m+2$ is slightly more complicated. Indeed, according to Theorem \ref{eq:deltas}, in order to satisfy condition \eqref{eq:gr4m+2}, we end up with $6$ linear homogeneous equations and only $4$ unknowns. But it can be proved that, among these $6$ relations, only $4$ of them are independent.
\end{itemize}
 Moreover, the coefficients involve the Weierstrass $\wp$ function at points of the type $k\omega_3$, $k$ integer, which can be computed via standard recursive schemes. The proof of the theorem is concluded.
 \fproof
 \subsubsection{About expressing the criterion in terms of the coefficients of the matrix $\Pb$}
By construction, we can a priori write 
\begin{equation}\label{eq:deltas1}
\delta^s(x) = M_s(x)y +N_s(x) \ \pmod{Q(x,y)},
\end{equation}
where $M_s$ and $N_s$ are rational fractions, whose  numerators and denominators are polynomials of unknown degrees, but with coefficients given in terms of  polynomials of the entries of $\Pb$. Moreover the decomposition \eqref{eq:deltas1} is unique, so that, comparing with \eqref{eq:deltas}, we have
\[
 \begin{cases}
\DD M_s(x) = \frac{4q^2\wp'(s\omega_3) a(x)}{W_s(x)}, \\[0.2cm]
\DD N_s(x) = \frac{V_s(x)}{W_s(x)},
\end{cases}
\]
where  $V_s,W_s$ are the second degree polynomials given by Theorem \ref{th:deltas}.
Moreover, by homogeneity, we can always rewrite
\[
M_s(x) = \frac{A_sa(x)}{F_s(x)},
\]
where   $A_s, K_s$ are real constants with 
\[
A_s = K_s [4q^2\wp'(s\omega_3)], \quad F_s(x) = K_sW_s(x).
\]
 In particular, by Corollary \ref{orders}, the group is of order $4s$ if and only if  $M_s\equiv0$, or equivalently  $A_s=0$, where now $A_s$ \emph{depends only on the entries of} $\Pb$ in a complicated polynomial form.

When the group is of order $4s+2$, we can exchange the role of $x$ and $y$ in the uniformisation \eqref{eq:xunif}, by uniformizing indeed $y(\omega)$. Then, \emph{mutatis mutandis}, this yields
\[
\delta^s(x)= \frac{\wt{A}_s\wt{a}(y)x + \wt{V}_s(y)}{\wt{F}(y)}\  \pmod{Q(x,y)},
\]
where $\wt{F}_s, \wt{V}_s$ are second degree polynomials. Referring again to Corollary \ref{orders}, we conclude that the group is of order $4s+2$ if and only if $\wt{A}_s=0$, where $\wt{A}_s$ depends on the coefficients of $\Pb$ in a polynomial way.

Finally we have shown that the finiteness of the the group is always equivalent to the cancellation of a \emph{single constant}, and this condition is tantamount to the existence of an algebraic hypersurface in the space of parameters. 

 \section{Some examples}

\subsection{$\H$ of order 4}
\begin{enumerate}
\item The product of 2 independent random walks inside the quarter
plane, so that
$$\sum_{i,j} p_{ij} x^i y^j = p(x) \wt{p}(y).$$
\item The {\em simple} random walk where $\sum p_{ij} x^i y^j = p(x) +
\wt{p}(y)$. Thus, in the interior of the quarter plane,
$p_{ij}\neq 0$ if, and only if, $i$ or $j$ is zero.
\item The case (a) in figure \ref{fig4.1.1}, which can be viewed as a
simple queueing network with parallel arrivals and internal transfers.
\end{enumerate} 

\subsection{$\H$ of order 6} 
\begin{figure}[htb]
\vspace{1cm}
\begin{center} 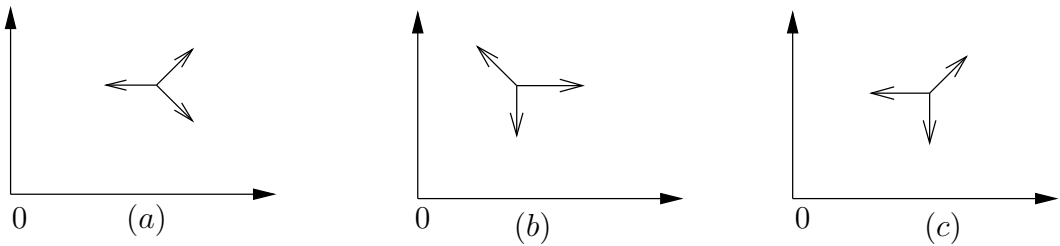 \end{center}
\caption{ \label{fig4.1.1} Some groups of order 6}
\end{figure}

Among these are the
cases studied in \cite{FAIA79}, \cite{BLAN82} and \cite{FLAT84}, which
are represented in figure \ref{fig4.1.1} (b), (c) respectively, where,
as before, only jumps inside the quarter plane have been drawn.

\appendix 
\section{About the uniformization of \protect $Q(x,y)=0$ (see \cite{FIM}, chapter 3)}  \label{app:unif}

Let $X(y)$ [resp.\ $Y(x)$] be the algebraic function defined by $Q(X(y),y)=0$ [resp.\ $Q(x,Y(x))=0$]. This function has two \emph{branches}, say $X_0$ and $X_1$ [resp.\ $Y_0$ and $Y_1$]. 

We shall denote by $\{x_\ell\}_{1\leq \ell\leq 4}$ the four roots of $D(x)$, which are the branch points of the Riemann surface 
	\begin{equation*}
          \mathscr{K}=\{(x,y)\in\mathbb{C}^{2} : Q(x,y) = 0\}.     
     \end{equation*}
 They are enumerated in such a way that $|x_1|\leq|x_2|\leq|x_3|\leq|x_4|$. 
 
 Moreover $x_1\leq x_2$, $[x_1x_2] \subset [-1,+1]$ and $0\leq x_2\leq x_3$.

When the associated Riemann surface is of genus $1$, the algebraic curve $Q(x,y)=0$ admits a \emph{uniformization} given in terms of the Weierstrass $\wp$ function with periods 
$\omega_1,\omega_2$ and its derivatives. Indeed, letting 
\begin{eqnarray*}
D(x) & = & b^2(x) -4a(x)c(x) \ \egaldef \ d_{4}x^{4}+d_{3}x^{3}+d_{2}x^{2}+d_{1}x+d_{0},\\
z & \egaldef & 2a(x)y+b(x),
\end{eqnarray*}
Then the following formulae hold (see \cite{FIM}, Section 3.3).
 \begin{enumerate}
\item If $d_4 \neq 0$ ($4$ finite branch points $x_1,\ldots x_4$) then $D'(x_4)>0$ and
\begin{equation}\label{eq:xunif}
\begin{cases}
x(\omega) = x_4 + \dfrac{D'(x_4)}{\wp(\omega) - \frac{1}{6} D''(x_4)},
\\[3ex] z(\omega) = \dfrac{D'(x_4) \wp'(\omega)}{2\left(\wp(\omega) -
\frac{1}{6} D''(x_4)\right)^2}. 
\end{cases}
\end{equation}
\item If $d_4 = 0$ ($3$ finite branch points $x_1,x_2,x_3$ and $x_4=\infty$) then
\begin{equation*}
\begin{cases}
x(\omega) = \dfrac{\wp(\omega) - \dfrac{d_2}{3}}{d_3}, \\[2.5ex]
z(\omega) = - \dfrac{\wp'(\omega)}{2d_3}. 
\end{cases}
\end{equation*}
\end{enumerate}
We also have
   \begin{equation} \label{eq_omega123}
          \omega_1 = 2i \int_{x_1}^{x_2}\frac{\text{d}x}{\sqrt{-D(x)}},
          \qquad \omega_2 = 2\int_{x_2}^{x_3}\frac{\text{d}x}{\sqrt{D(x)}},
          \qquad \omega_3 = 2\int_{X(y_1)}^{x_1}\frac{\text{d}x}{\sqrt{D(x)}},
     \end{equation}
noting that $\omega_1$ is purely imaginary, while $0<\omega_3<\omega_2$.  It was proved  that  the group $\H$ is finite of order $2n$ if, and only if,
\[
 n \omega_3 = 0 \quad \bmod (\omega_1,\omega_2), 
 \] 
 or, since
$\omega_3$ is real,
\begin{equation} \label{eq:gfinite}
n \omega_3 = 0 \quad \bmod (\omega_2), 
\end{equation} where
$n$ stands for the minimal positive integer with this property. 
On the \emph{universal covering} $\C$ (the finite complex plane), the automorphisms introduced in Section \ref{GENUS} become (see \cite{FIM}, Section 3.3)
  \begin{equation} \label{eq:auto}
          \xi^*(\omega) = -\omega+\omega_2,
          \qquad \eta^*(\omega) = -\omega + \omega_2 + \omega_3,
          \qquad \delta^*(\omega) = \eta^*\xi^*= \omega+\omega_3,
     \end{equation}
and here $\delta=\xi\eta$ corrresponds to $\delta^*= \eta^*\xi^*$ (note the permutation of the letters\,!). For any $f(x,y) \in C_Q(x,y)$, 
\[
\delta(f(x,y)) = f(\delta(x),\delta(y) = f(x(\delta^*(\omega)), y(\delta^*(\omega))), \ \omega\in\C.
\]
In particular,
\begin{equation}\label{eq:auto1}
\begin{cases}
\xi\eta(x)  =  x(\eta^*\xi^*(\omega))= x(\omega+\omega_3), \\[0.1cm]
\eta(x)  =  x(\eta^*(\omega))= x(-\omega+\omega_2+\omega_3) = x(\omega-\omega_3).
\end{cases} 
\end{equation}   

\section{Some symmetric quantities of the $\wp$ function}\label{app:calculs}
From the form of $x(\omega)$ in \eqref{eq:xunif} and the expressions \eqref{eq:auto} of the automorphisms $\xi^*,\eta^*,\delta^*$, it appears  we need to calculate the quantities 
\[
A(u,v) \egaldef \wp(u+v) + \wp(u-v), \quad B(u,v) \egaldef  \wp(u+v)\wp(u-v)
\]
in terms of rational functions of  $\wp(u)$ and $\wp(v)$. 

Letting $\wp'(u)$ denote the derivative of $\wp(u)$ with respect to $u$,  it is well known (see e.g. \cite{JS}) that  $\wp$ is even, $\wp'(z)$ is odd, and they satisfy the algebraic differential equation
\begin{equation}\label{eq:wp}
 \wp'^2 = 4\wp^3 -g_2\wp - g_3,
 \end{equation}
 where the quantities $g_2, g_3$ are constants often called \emph{the invariants}. Equation \eqref{eq:wp} admits also of  the classical form
 \begin{equation}\label{eq:wp1}
  \wp'^2 = 4(\wp-e_1)(\wp-e_2)(\wp-e_3),
  \end{equation}
  where 
  \[
  e_1+e_2+e_3=0, \ g_2= -4(e_1\,e_2 + e_2\,e_3 + e_3\,e_1), \ g_3 = 4e_1\,e_2\,e_3,
  \]
\[
e_1= \wp\left(\frac{\omega_2}{2}\right), \quad 
e_2 = \wp\left(\frac{\omega_1+\omega_2}{2}\right), \quad e_3= \wp\left(\frac{\omega_1}{2}\right).
\]
In the present situation, all the $e_i$'s are real and satisfy 
 \begin{equation}\label{eq:wpe}
 e_1>0, \quad  e_3<0, \quad e_1> e_2 > e_3.
 \end{equation}

We shall also need the \emph{addition theorem} (see \cite{JS}), which reads, for all $u,v, u\pm v\neq0$,
\[
\wp(u+v) = \frac{1}{4} \biggl(\frac{\wp'(u)-\wp'(v)}{\wp(u)-\wp(v)}\biggr)^2 -\wp(u)-\wp(v).
\]
It is important to note that $A(u,v)$ and $B(u,v)$ are symmetric and even w.r.t. $(u,v)$.

Let us set for a while $X\egaldef \wp(u)$,  $Y\egaldef \wp(v)$. 
\begin{lem}\label{lem:WSP}
\begin{eqnarray}
A(u,v) &=& \frac{(X + Y)(4XY -g_2)-2g_3}{2(X-Y)^2}, \label{eq:swp}\\[0.2cm]
B(u,v) &=& \frac{(XY)^2 + \frac{g_2}{2}XY + g_3(X+Y) + \frac{g_2^2}{16}}{(X-Y)^2}. \label{eq:pwp}
\end{eqnarray}
\end{lem}
\dproof
The parity properties allow to write
\[
A(u,v) = \frac{\wp'^2(u) + \wp'^2(v) - 4(\wp^2(u) -\wp^2(v))(\wp(u)-\wp(v))}{2(\wp(u)-\wp(v))^2},
\]
which by \eqref{eq:wp} yields directly \eqref{eq:swp}.
As for $B(u,v)$, we have
\begin{equation}\label{eq:p1wp}
B(u,v)= \frac { \left[(X+Y)(4XY - g_2) - 2g_3 \right] ^{2}-4\, (4 X^3 - g_2 X-g_3) (4 Y^3 - g_2 Y-g_3)} 
{16\, (X - Y) ^4}.
\end{equation}
Then it is quite reasonable to guess the last expression reduces in fact to a  rational fraction of second degree. This is clear in our particular context, where we shall take for example $u=\omega$ and 
$v=m\omega_3/2$, but it turns out to be true for arbitrary $u,v$. Indeed, the numerator of \eqref{eq:p1wp} can be factorized by  $(X-Y)^2$, which leads to \eqref{eq:pwp} (this last property has been enlighted with the help of the  \emph{Maple~18} Computer Algebra System).

The proof of the lemma is concluded.
\fproof

Now, by Lemma \ref{lem:WSP} and equation \eqref{eq:xunif},  we are in a position to compute 
\begin{equation}\label{eq:SP}
S(u,v) \egaldef x(u+v)+x(u-v), \quad  P(u,v) \egaldef x(u+v)x(u-v).
\end{equation}
For the sake of brevity, It will be convenient to write 
\begin{equation}\label{eq:xunif1}
x(\omega) = p + \dfrac{q}{\wp(\omega) - r},
\end{equation}
 where $p,q,r$  are known constants coming in equation \eqref{eq:xunif}.
Then  the following functional algebraic relations, giving $S(u,v)$ and $P(u,v)$ in terms of $A(u,v)$ and $B(u,v)$, are straightforward. 
\begin{equation}\label{eq:SP1}
\begin{cases}
\DD S  = \frac{2pB + (q-2pr)A + 2r(pr-q)}{B-rA+r^2},\\[0.3cm]
\DD P = \frac{p^2B + p(q-pr)A + (pr-q)^2}{B-rA+r^2}.
\end{cases}
\end{equation}
Ad libitum, we shall specify  the variables only whenever needed: for instance, instead of $S$, we shall write $S(X,Y)$, or $S(\omega,s\omega_3)$, etc. So, by using  \eqref{eq:swp}, \eqref{eq:pwp}, \eqref{eq:SP1}, we obtain the final expressions of $S$ and  $P$, which, as expected, are rational functions of \emph{second degree} with respect to $X$ and $Y$ separately, and thus also with respect to $x(u)$ and $x(v)$.  

Moreover, setting 
\begin{equation}\label{eq:SP2}
A_1= 2(X-Y)^2A, \quad B_1=2(X-Y)^2B, \quad D=2(X-Y)^2(B-rA+r^2),
\end{equation}
 one sees immediately that $S$ and $P$ can be expressed as ratios of polynomials of second degree  in $x$ through the homography \eqref{eq:xunif1}.

\paragraph{Acknowledgements}
The authors would like to thank Frederic Chyzak (Inria-Saclay) for his valuable help concerning the use of the function \emph{RootOf}  in the \emph{Maple~18} Computer Algebra System.


\end{document}

%% file: order6.pdf_t
\begin{picture}(0,0)%
\includegraphics{order6}%
\end{picture}%
\setlength{\unitlength}{4144sp}%
\begingroup\makeatletter\ifx\SetFigFont\undefined%
\gdef\SetFigFont#1#2#3#4#5{%
  \reset@font\fontsize{#1}{#2pt}%
  \fontfamily{#3}\fontseries{#4}\fontshape{#5}%
  \selectfont}%
\fi\endgroup%
\begin{picture}(6260,1460)(670,-1035)
\put(3106,-916){\makebox(0,0)[lb]{\smash{{\SetFigFont{12}{14.4}{\familydefault}{\mddefault}{\updefault}0}}}}
\put(6121,-961){\makebox(0,0)[lb]{\smash{{\SetFigFont{12}{14.4}{\familydefault}{\mddefault}{\updefault}$(c)$}}}}
\put(1397,-917){\makebox(0,0)[lb]{\smash{{\SetFigFont{12}{14.4}{\familydefault}{\mddefault}{\updefault}$(a)$}}}}
\put(3692,-962){\makebox(0,0)[lb]{\smash{{\SetFigFont{12}{14.4}{\familydefault}{\mddefault}{\updefault}$(b)$}}}}
\put(721,-916){\makebox(0,0)[lb]{\smash{{\SetFigFont{12}{14.4}{\familydefault}{\mddefault}{\updefault}0}}}}
\put(5356,-916){\makebox(0,0)[lb]{\smash{{\SetFigFont{12}{14.4}{\familydefault}{\mddefault}{\updefault}0}}}}
\end{picture}%